\documentclass{article}
\usepackage{amsmath, amssymb, a4, tabbert}

\newcommand {\supp}{\text{supp}\,}
\begin{document}
\bibliographystyle{alpha}
\title{A Remark on the Leray's Problem on Stationary
Navier-Stokes Flows with Large Fluxes in Infinite Cylindrical
Domains}
\author{Myong-Hwan Ri$^\dag$}
\date{}
\maketitle
\begin{abstract}
\noindent
We consider Leray's problem on stationary Navier-Stokes flows with arbitrary large fluxes in
an unbounded cylinder with several exits to infinity.
For a stationary Navier-Stokes flow with
large fluxes in the unbounded cylinder in the sense of Definition \ref{D2.1},
we prove that, if
the  difference between the pressure of the main flow and
the pressure of the Poiseuille flow with the same flux
in a branch of the cylinder
satisfies some asymptotic boundedness condition at $|x|\ra\infty$, see \eq{4.0},
then the flow behaves at infinity  of the branch like the Poiseuille flow.

\end{abstract}
\vpn{\small{\bf 2000 Mathematical Subject Classification:} 35Q30;
35B35;76D05; 76D07; 76E99 
\let\thefootnote\relax\footnote{\hspace{-0.3cm}
\dag: Institute of Mathematics, State Academy of Sciences, DPR Korea.
He was supported by 2012 CAS-TWAS Postdoctoral Fellowship with grant number 3240267229}.

\section {Introduction and main result}

Let
\begin{equation}
\label{E1.1} \Om=\bigcup_{i=0}^{m}\Om^i
\end{equation}
be a cylindrical domain of $C^{2}$-class of $\R^3$, where $\Om^0$
is a bounded domain and $\Om^i, i=1,\ldots, m,$ are disjoint
semi-infinite straight cylinders, that is, in possibly different
coordinates,
$$ \Om^i=\{x^i=(x^i_1,x^i_2, x^i_3)\in\R^3: \:x^i_3>0,
   x'^{i}=(x^i_1,x^i_{2})\in\Si^i\},$$
where $\Si^i\subset \R^{2}$, $i=1,\ldots,m,$ is a bounded domain and
$\Om^i\cap\Om^j=\emptyset$ for $i\neq j$. Without loss of
generality, we assume for each $i=1,\ldots, m$ that the coordinate
system which is fixed in $\Om^i$ is such that $x'^i, x_3^i$ denote
the variables with respect to the cross section $\Si^i$ and the
axial direction of $\Om^i$, respectively.

Let us consider the stationary Navier-Stokes system
\begin{equation}
\label{E1.2}
\begin{array}{rclll}
     -\Da U  + (U\cdot\na)U+\na P & =& 0 \,\, &\text{in}&\Om,\ek
     \div U &=& 0 \,\, &
   \text{in} &\Om,\ek
     U &=& 0 \,\, &\text{on}&\pa\Om.
\end{array}
\end{equation}
We impose an additional condition for the behavior of the velocity
field at infinity as
 \begin{equation}
 \label{E1.2n}
 \di\lim_{|x|\ra\infty}U(x)=u_\infty,
 \end{equation}
where $u_\infty$ coincides at infinity of each exit
$\Om^i,i=1,\ldots,m,$ with the {\it Poiseuille flow} ${\bf v}_i$ in
$\Om^i$ corresponding to the prescribed flux $\Phi_i$.
\par Poiseuille flows in an infinite
straight cylinder $\Si\ti\R$ are often referred to as flows parallel
to the axial direction. In the stationary case, the Poiseuille flow
${\bf v}=(0,0,v(x'))$ and the corresponding {\it Poiseuille flow pressure}
$\Pi(x_3)=-kx_3+b$ are simply given by the Poisson equation
$$-\Da'  v=k,\quad  v|_{\pa \Si}=0,$$
where, if $\Si$ is a Lipschitz domain, then $k=c(\Si)\Phi$,
 where $c(\Si)=\frac{1}{\int_\Si |\na' g|^2\,dx'}$,
$-\Da'g=1,\; g|_{\pa\Si}=0$.

Note, due to
the solenoidal condition for the fluid, that if $U$ satisfies \eq{1.2n},
then
\begin{equation}
\label{E1.4}
 \int_{\Si^i}U \cdot {\bf n}^i\,dx'^i=\Phi_i
 \end{equation}
should necessarily hold true, where ${\bf n}^i$ is the unit vector
along the positive axial direction of $\Om^i$.  Moreover, the flux
$\Phi_i$ should be independent of $x_3^i$ over $\Om^i$ for
$i=1,\ldots,m$ and
\begin{equation}
\label{E1.5}
 \sum_{i=1}^m \Phi_i=0
\end{equation}
should be naturally assumed.

Classical {\it Leray's problem} is to show whether or not the problem \eq{1.2}-\eq{1.4}
will admit a solution. Leray's problem seems to have been proposed, see \cite{Am77},
 by J. Leray himself to O.A. Ladyzhenskaya,
who in \cite{La59} attempted an existence proof under no restrictions on
the viscosity.

There is a number of
papers dealing with stationary Leray's problem.
Fundamental contribution to Leray's problem was made by Amick in \cite{Am77},
where the existence of unique weak solution to \eq{1.2}-\eq{1.4} was proved
under a  smallness assumption on the total flux $\sum_{i=1}^m|\Phi_i|$,
see also \cite{AmFr80}, \cite{Fa03}, \cite{FoFr00}, \cite{Ka85}, \cite{KaPi84},
and \cite{NaPi84}-\cite{RF07}.
However, it has been shown, up to now, that Leray's problem is solved positively
only under smallness assumptions on the total flux, and the problem
for arbitrary
large total flux is known as one of the most challenging problems in the
theoretical fluid dynamics;
for the Leray’s and related problems we refer, in particular, to
\cite{Ga94-1}, Chap. VI, Sections 1 and 2, and \cite{Ga94-2}, Chap. XI, Sections 1, 2, 3 and 4,
cf. also \cite{Be04},
Introduction and references cited therein for more details.
\par\bigskip
In this paper, we aim at considering the Leray's problem for large total flux;
we present a condition on the flow pressure to allow a weak solution \eq{1.2}, \eq{1.4}
to behave like Poiseuille flows at $|x|\ra\infty$.

In order to explain the main result of the paper,
let us give the definition of the weak solution
to the system \eq{1.2}, \eq{1.4}.
Let $$C^\infty_{0,\si}(\Om):= \{\vp\in C^\infty_0(\Om)^n: \,\,\div\vp=0\}$$
 and domains $\Om^i_N$
and $\Om_{N,N+1}$ be respectively given as
$$\Om^i_N=\{x\in\Om^i:x^i_3\leq N\},\,\,
  \Om^i_{N,N+1}=\{x\in\Om^i:N \leq |x|\leq N+1\}, \,\, \Om_{N,N+1}=\cup_{i=1}^m\Om^i_{N,N+1}.$$

\begin{tdefi}
\label{D2.1} {\rm A vector field $U:\Om\ra \R^3$
is called a weak solution to \eq{1.2},\eq{1.4} if it
satisfies the conditions (i)$\sim$(v):
\begin{itemize}
 \item[(i)] $U\in H^1_{\text{loc}}(\bar\Om)$ satisfies
 \begin{equation}
\label{E1.6}
\begin{array}{l}
\di\int_{\Om^i_N}|\na U|^2\,dx\leq K_1N \quad\text{ for all } N>0\\[3ex]
\di\int_{\Om_{N,N+1}}|\na U|^2\,dx\leq K_2 \quad\text{ for all } N>0
\end{array}
\end{equation}
with $K_1$ and $K_2$ independent of $N$.
 \item[(ii)]
The variational equation
\begin{equation}
\label{E2.1} (\na U,\na \psi)=(U\cdot\na\psi, U)\quad\text{for all }\psi\in C^\infty_{0,\si}(\Om)
\end{equation}
holds.
\item[(iii)] $U$ vanishes on the boundary $\pa\Om$.
 \item[(iv)] Solenoidal condition $\div U = 0$ in $\Om$ holds in the distributional sense.
\item[(v)] $U$ satisfies \eq{1.4} in the trace sense.
 \end{itemize}
 }
\end{tdefi}

In \cite{LS83}
the existence of a weak solutions to \eq{1.2}, \eq{1.4} in the sense of
Definition \ref{D2.1} was proved without any smallness assumption on the total flux.
Note that if $U$ is a solution to \eq{1.2}, \eq{1.4} in the sense of Definition \ref{D2.1}, then
there is an associated pressure $P$, which is determined
uniquely up to a constant difference, such that $(U,P)$ solves the system \eq{1.2}
in the sense of distribution. Hence, we shall also call $(U,P)$ a weak solution
to \eq{1.2}, \eq{1.4}.

Even if a weak solution $U$ to \eq{1.2}, \eq{1.4} in the sense of
Definition \ref{D2.1} exists for any large flux,
it is not known yet whether the flow $U$ will tend
to Poiseuille flows corresponding to given fluxes in each exit of $\Om$
as $|x|\ra\infty$.
\par\bigskip
The main result of this paper is the following statement:


\begin{theo}
\label{T1.2}
{\rm
Let $\Pi_i$ be pressure of
Poiseuille flow with flux $\Phi_i$ in $\Om^i$ for $i=1,\ldots,m$
satisfying the flux condition \eq{1.5}.
Let $(U,P)$ be any weak solution to the problem \eq{1.2}, \eq{1.4}
in the sense of Definition \ref{D2.1}.

If
\begin{equation}
\label{E4.0}
\liminf_{N\ra\infty}\int_{\Si^i} |(P-\Pi_i)(x'^i,N)|\,dx'^i<\infty, i=1,\ldots,m,
\end{equation}
then $U-{\bf v}_i\in H^1(\Om_i), i=1,\ldots,m,$ and, in particular,
the stationary solution $U$ behaves at $|x|\ra\infty$ like
Poiseuille flow ${\bf v}_i$  in each branch $\Om_i,i=1,\ldots,m$, of $\Om$.
}
\end{theo}
\begin{rem}
Obviously, if $U$ behaves like a Poiseuille flow ${\bf v}_i$ at $|x|\ra \infty$ in $i$-th exit of $\Om$,
then the pressure $P$ of the flow also behaves
like the corresponding pressure $\Pi_i$ of the Poiseuille flow and hence
the relation
$$\lim_{|N|\ra \infty}\int_{\Si^i}(P-\Pi_i)(x'^i,N)\,dx'^i=0$$
must hold true.
\end{rem}

\begin{rem}
From Theorem \ref{T1.2} one can get the following deduction:
Suppose a fluid flow $U=(U_1,U_2,U_3)$ (with pressure $P$ and)
with nonzero flux $\Phi_i$ in each $\Om_i$ satisfies the assumptions of Definition \ref{D2.1}.
Suppose $U$ does not behave at $|x|\ra \infty$ like Poiseuille flow ${\bf v}_i=(0,0,v_i)$ (denote its pressure by $\Pi_i$) with flux $\Phi_i$ in $\Om_i$ for some $i\in \{1,\ldots,m\}$.
It is seen that \eq{4.0} is equivalent to
\begin{equation}
\label{E1.9}
\liminf_{N\ra\infty} \Big|\int_1^N \int_{\pa\Si^i} \frac{\pa (U_3-v_i)}{\pa n'}\,dx'\Big|<\infty,
\end{equation}
where $n'$ is the unit outward normal vector at $\pa\Si^i$,
see \eq{4.24}.
We notice that
$\int_1^N \int_{\pa\Si} \frac{\pa (U_3-v)}{\pa n'}\,dx'$ is the
axial directional component of friction force caused in the part $\Om^i_{1,N}$
by the velocity perturbation $U-{\bf v}_i$.

Therefore, one may conclude from Theorem \ref{T1.2} that
if the axial directional component of
frictional force (resistance in the flow direction)
caused by the perturbation flow from the Poiseuille flow in $\Om_i$ is finite, then
the flow $U$ with uniform property as in Definition \ref{D2.1}
behaves at $|x|\ra\infty$ like the Poiseuille flow ${\bf v}_i$ in $\Om_i$.
\end{rem}

We use the following notations.

For a domain $G$ of $\R^k,k\in \N,$
let $L^r(G), W^{s,r}(G)$, $s>0$, $1<r \leq \infty,$ be the usual Lebesgue and Sobolev
spaces on $G$.
The space  $H^1_{0}(G)$
($W^{1,2}_0(G)$) is the $H^{1}$-completion of the set $C^\infty_{0}(G)$.

 As long as no confusion arises, we use the same notations for scalar and vector spaces
and for constants, e.g. $c, C$, appearing in the proofs.

%
%
%
\section{Preliminaries on the weak solution}

In this section we show some properties of the weak solution $(U,P)$ to \eq{1.2},\eq{1.4}.
Let $\Phi=\sum_{j=1}^m|\Phi_j|$.

\begin{propo}
\label{P3.1} {\rm
   Let $(U,P)$ be a weak solution to \eq{1.2}, \eq{1.4} in the sense of Definition \ref{D2.1}.
Then, there holds
$$U, \na^2U,\na P\in L^r(\Om_{N,N+1}), \,\,
\|U, \na^2U, \na P\|_{L^r(\Om_{N,N+1})}\leq C(r,\Om,\Phi), \,\,
\forall r\in(1,\infty),\,\,
\forall N\in \N,$$
with a constant $C>0$ independent of $N\in\N$.}
\end{propo}
{\bf Proof}:
 Fix an arbitrary index $i\in \{1,\ldots,m\}$.
 Without the loss of generality
 we may assume that the coordinate system $x$ coincides with the one $x^i$ fixed in $\Om^i$.
 Given $N\in\N\cup\{0\}$ let
$$G_N:=\Om^i_N \quad \text{for }N\in\N,\quad G_0: = \emptyset,\quad G_{N,N+k}:=G_{N+k}\setminus G_{N}, k\in\N.$$
Note that by \eq{1.6} one has
\begin{equation}
\label{E3.5} U|_{G_{N,N+1}}\in H^1(G_{N,N+1}),\quad
\|U\|_{H^1(G_{N,N+1})} \leq K_2,\quad \forall N\in \N.
\end{equation}
Now, fix a domain $G$ of $C^2$-class such that
$$G_{1,2}\subset G\subset G_{0,3}$$
and $\pa G\cap\pa \Om$ is a strictly interior of $\pa G_{0,3}\cap \pa\Om$.
Then,
$$G_{N,N+1}\subset G^{(N)}\subset G_{N-1,N+2},$$
where $G^{(N)}$ is obtained by shifting the domain $G$ as the distance $N$ in the positive direction of the
axis of $\Om^i$, that is,
$$G^{(N)}:=G+\{(0,0,N)\}.$$
If the boundary of $\Om$ is smooth enough, then one may get
$U, P\in C^2(\bar G^{(N)})$
by  Chap. XI, Theorem 1.1 in \cite{Ga94-2}. However, since we do not have
this smoothness for the boundary, we need a more refined argument
for the estimate of $(U, P)$.

In view of the geometry of $G_{N-1,N+2}$, we get by Sobolev embedding theorem and \eq{3.5}
that
\begin{equation}
\label{E3.6} \begin{array}{l}
 (U\cdot\na)U\in
L^{3/2}(G_{N-1,N+2}),\ek \|(U\cdot\na)U\|_{L^{3/2}(G_{N-1,N+2})}\leq
c\|U\|^2_{H^1(G_{N-1,N+2})}\leq c_1(\Si^i,\Om,\Phi)
\end{array}
\end{equation}
with constants $c, c_1$ depending only on $\Si^i$
and independent of $N$.
Moreover, we have
\begin{equation}
\begin{array}{rcl}
\label{E3.7} \|\na P\|_{H^{-1}(G_{N-1,N+2})} &\leq & \|\Da
U-(U\cdot\na)U\|_{H^{-1}(G_{N-1,N+2})}\ek
        &\leq & c(\|\na U\|_{L^2(G_{N-1,N+2})}+
        \|(U\cdot\na)U\|_{L^{3/2}(G_{N-1,N+2})})\ek
        &\leq &c(\Si^i,\Om,\Phi),
\end{array}
\end{equation}
where the constants are independent of $N\in \N$, since $H^1_0(G_ {N-1,N+2})$ is continuously
embedded into $L^3(G_{N-1,N+2})$ with embedding constant depending only on $\Si^i$
and independent of $N$.
\par
 Let $P_N=\frac{1}{|G_{N-1,N+2}|}\int_{G_{N-1,N+2}}P\,dx$.
Since $(U,P-P_N)$ solves the system
\begin{equation}
\label{E3.8}
\begin{array}{rccl}
  -\Da U +\na (P-P_N) &=&-(U\cdot\na)U \,\,
  &\text{in }G,\ek
  \div U&=&0 \,\, &\text{in }G,\ek
  U&=&0\,\,&\text{on }\pa G_{N-1,N+2}\cap \pa \Om,
  \end{array}
  \end{equation}
\cite{Ga94-1}, Chap. IV, Theorem 5.1 implies
$U\in W^{2,3/2}(G^{(N)}), P-P_N\in W^{1,3/2}(G^{(N)})$
and
\begin{equation}
\label{E3.9}
\begin{array}{l}
\|U\|_{W^{2,3/2}(G^{(N)})}+\|P-P_N\|_{W^{1,3/2}(G^{(N)})}\ek \leq
C_1\big(\|(U\cdot\na)U\|_{L^{3/2}(G_{N-1,N+2})}+
\|U\|_{W^{1,3/2}(G_{N-1,N+2})}+\|P-P_N\|_{L^{3/2}(G_{N-1,N+2})}\big)\ek
\leq C_1C_2\big(\|(U\cdot\na)U\|_{L^{3/2}(G_{N-1,N+2})}+
\|U\|_{H^{1}(G_{N-1,N+2})}+\|P-P_N\|_{L^{2}(G_{N-1,N+2})}\big)
\end{array}
\end{equation}
with  $C_1$ and $C_2$ independent of $N\in\N$.
In fact, by \cite{Ga94-1}, Chap. IV, Theorem 5.1,
$C_1$ depends only on $G^{(N)},G_{N-1,N+2}$ and hence
$$C_1=C_1(G,G_{0,3}),\quad
C_2=\max\{1,|G_{N-1,N+2}|^{1/3}\}=\max\{1,(3|\Si^i|)^{1/3}\}$$
 since $G^{(N)}$ and
$G_{N-1,N+2}$ are  obtained by shifting $G$ and $G_{0,3}$, respectively.
\par
On the other hand, since $P-P_N$ has mean value $0$ in $G_{N-1,N+2}$,
we get by \cite{Bog79}, Theorem 1 that
\begin{equation}
\label{E3.10} \|P-P_N\|_{L^{2}(G_{N-1,N+2})}\leq c\|\na
P\|_{H^{-1}(G_{N-1,N+2})},
\end{equation}
where the constant $c$ depends on the diameter of $G_{N-1,N+2}$
and consequently does not depend on $N$.

Thus, \eq{3.6}, \eq{3.7}, \eq{3.9} and \eq{3.10} yield
\begin{equation}
\label{E3.11}
\|U\|_{W^{2,3/2}(G^{(N)})}+\|P-P_N\|_{W^{1,3/2}(G^{(N)})}\leq
C(\Si^i,\Om, \Phi).
\end{equation}
\par
Note that \eq{3.11} holds for the sets $G^{N-1}$ and $G^{N+1}$ as well.
Therefore, we have
\begin{equation}
\label{E3.11n}
\|U\|_{W^{2,3/2}(G_{N-1,N+2})}+\|P-P_N\|_{W^{1,3/2}(G_{N-1,N+2})}\leq
C(\Si^i,\Om, \Phi).
\end{equation}
\par
Repeating the above argument, we get $(U\cdot\na)U\in L^{q}(G_{N-1,N+2})$ and
$$
\|(U\cdot\na)U, P-P_N\|_{L^q(G_{N-1,N+2})} \leq
c(q,\Si^i,\Om,\Phi),\quad\forall q\in (1,3),$$
with constant $c>0$ independent of $N$.
Hence $U\in W^{2,q}(G^{(N)})$,
$P-P_N\in W^{1,q}(G^{(N)})$ and we have
\begin{equation}
\label{}
\begin{array}{l}
\|U\|_{W^{2,q}(G^{(N)})}+\|P-P_N\|_{W^{1,q}(G^{(N)})}\ek \leq
C_3\big(\|(U\cdot\na)U\|_{L^{q}(G_{N-1,N+2})}+
\|U\|_{W^{1,q}(G_{N-1,N+2})}+\|P-P_N\|_{L^{q}(G_{N-1,N+2})}\big)\ek
\leq C_4,
\end{array}
\end{equation}
where $C_3,C_4$ depend on $q,G_{0,3},G,\Si^i,\Om,\Phi$.
Applying  the above argument once again,
in view of the continuous embedding  $W^{1,q(r)}\hookrightarrow L^r$  ($\forall r\in
(1,\infty)$  $\exists q(r)\in (1,3)$), we get finally that
 $U\in W^{2,r}(G^{(N)}), P-P_N\in
W^{1,r}(G^{(N)}),\forall r\in (1,\infty)$
and
\begin{equation}
\label{}
\begin{array}{l}
\|U\|_{W^{2,r}(G^{(N)})}+\|P-P_N\|_{W^{1,r}(G^{(N)})}\ek \leq
C_5\big(\|(U\cdot\na)U\|_{L^{r}(G_{N-1,N+2})}+
\|U\|_{W^{1,r}(G_{N-1,N+2})}+\|P-P_N\|_{L^{r}(G_{N-1,N+2})}\big)\ek
\leq C_6,
\end{array}
\end{equation}
where $C_5, C_6$ depend only on $q,G_{0,3},G,\Si^i,\Om,\Phi$.

The proof of the proposition is complete. \qed

\begin{coro}
\label{C3.2} {\rm For $(U,P)$ let us assume the same as in Proposition \ref{P3.1}. Then,
$$\|U,\na U\|_{L^\infty(\Om)}\leq
C(\Om,\Phi).$$}
\end{coro}

\section{\large Poiseuille flow-like behavior of fluid flow at $|x|\ra\infty$}

Let us introduce cut-off functions
$\{\varphi_i\}_{i=0}^{m}$ for
 $\Om$ such that
\begin{equation}
\label{}
\begin{array}{l}
    \sum_{i=0}^{m}\varphi_i(x)=1, \quad 0\leq \varphi_i(x)\leq 1\quad
    \text{for }x\in\Om,\ek
\varphi_i \in C^{\infty}(\bar{\Om}^i), \quad
   \text{dist}\,(\text{supp}\,\varphi_i,\, \pa\Om_i\cap \Om)\geq \delta>0,
   \,\,i=0,\ldots, m,
\end{array}
\end{equation}
where 'dist' means the distance.
For $i=1,\ldots,m$ let ${\bf v}_i=(0,0,v_j)$ be the Poiseuille flow
with flux $\Phi_i$ in $\Si^i$.
We know that the corresponding pressure is given by
$\Pi_i=k_ix_3^i+b_i$, see Introduction.
Moreover, let $\tilde{\bf v}_i, \tilde\Pi_i$ be respectively
the zero extension of ${\bf v}_i, \Pi_i$ onto $\Om$.

Let ${\bf a}$ be a carrier of the Poiseuille flows ${\bf v}_i, i=1,\ldots,m,$
such that for all $r\in (1,\infty)$
\begin{equation}
\label{E3.1n} {\bf a}\in W^{2,r}_{\text{loc}}(\overline\Om),\,\,
\div\, {\bf a}=0 \;\text{in}\; \Om,\,\, {\bf a}=0\;
\text{on}\;\pa\Om,\,\, {\bf a}={\bf v}_i\;\text{in}\;\Om_i\setminus
\Om^0, i=1,\ldots,m.
\end{equation}
In \cite{Pi97}, see also \cite{Ga94-1}, Chap. 6, \S\S1, such a vector
field ${\bf a}$ was constructed as
$${\bf a=z+v},\quad {\bf v}=\sum_{i=1}^m\vp_i\s{\bf v}_i,
\quad {\bf z}\in W^{2,r}_0(\Om^0),$$
 Then,
\begin{equation}
\label{E3.1}
\|{\bf a}\|_{W^{2,r}(\Om^0)}+\|{\bf a}, \na {\bf
a}\|_{L^\infty(\Om)} \leq C(\Om)\Phi, \quad \forall r\in (1,\infty).
\end{equation}

Now, let
\begin{equation}
\label{E3.3n} u:=U-{\bf a},\quad p:=P-\tilde{P},
\quad \tilde{P}:=\sum_{i=1}^m\varphi_i\tilde{\Pi}_i.
\end{equation}
Then, $(u, p)$ solves
\begin{equation}
\label{E3.3}
\begin{array}{rcll}
  -\Da u  + (u\cdot\na){\bf a}+ ({\bf a}\cdot\na)u+(u\cdot\na)u+\na p &= &f \,\,
  &\text{in }\Om,\ek
  \div u &=&0 \,\, &\text{in }\Om,\ek
  u &=&0 \,\, &\text{on }\pa\Om,\ek
  \end{array}
\end{equation}
where $f=\Da {\bf a}-({\bf a}\cdot\na){\bf a}-\na \tilde{P}$.
From the construction of the vector field ${\bf a}$ and $\tilde{P}$ it is clear that
\begin{equation}
\label{E3.30}
 \supp f\subset \Om^0.
\end{equation}

%
{\bf Proof of the Theorem \ref{T1.2}:}
Throughout the proof, we fix $i\in \{1,\ldots,m\}$ and use
the same notation as in the proposition \ref{P3.1} for domains.
Let us assume w.l.o.g that
the coordinates $x^i$ coincides with the generic coordinate $x$
and $\Om^0 \cap G_{1,2} = \emptyset$.

 For each $N>0$ let us introduce a cut-off function $\eta=\eta_N$
such that
$$\eta \in C^\infty(\bar\Om),\quad \supp \eta\subset
\bar{G}_{N+1},\quad 0\leq \eta(x)\leq 1,\quad \forall x\in \Om,$$
 and for $x\in \Om^i$
$$\eta(x)=\eta(x_3)=\left\{\begin{array}{ll}
1& \text{ if }x_3\in(1, N)\ek 0& \text{ if }x_3\in (0,\ha)\cup
(N+\ha,\infty).
\end{array}\right.$$
Obviously, we have
$$\eta {\bf a}\in W^{2,r}(\Om)
\cap W^{1,r}_0(\Om^i_{N+1}),\quad \forall r\in (1,\infty).$$

Let
\begin{equation}
\label{E4.22}
u^N:=\eta u,
\end{equation}
then, obviously, $u^N\in H^1_0(\Om)$.
In the following we shall show that
\begin{equation}
\label{E4.2}
    \|u^{N}\|_{H^1_0(\Om)} \leq C(\Om,\Phi),\quad\forall
    N\in\N.
\end{equation}

   Testing \eq{3.3} with $\eta{\bf a}$ yields
\begin{equation}
\label{E4.1} (\na u,\na (\eta{\bf a}))+((u\cdot\na){\bf a},\eta{\bf
a}) +(({\bf a}\cdot\na)u,\eta{\bf a})+((u\cdot\na)u,\eta{\bf a})
-(p,\div(\eta{\bf a}))=(f,\eta{\bf a}),
\end{equation}
where and unless indicated $(\cdot,\cdot)$ stands for the scalar product in $G_{N+1}$.
By elementary calculation we have
\begin{equation}
\label{E4.3}
\begin{array}{l}
(\na u,\na (\eta{\bf a}))
 = (\na u^N,\na {\bf a})-(\na\eta u,\na{\bf a})
    +(\na u, \na\eta {\bf a})\ek
  =(\na u^N,\na {\bf a})+(\na u, \na\eta {\bf a})_{G_1\cup G_{N,N+1}}-(\pa_3\eta u, \pa_3 {\bf a})_{G_1}\ek
  =(\na' u_3,\na' v_i)_{G_{1,N}}+(\na' u^N_3,\na' v_i)_{G_1\cup G_{N,N+1}}
       +(\pa_3 u_3, \pa_3\eta v_i)_{G_1\cup G_{N,N+1}}\ek
       \hspace{0.5cm} -(\pa_3\eta u, \pa_3 {\bf a})_{G_1}\ek
  =(\na' u_3,\eta\na' v_i)_{G_1\cup G_{N,N+1}}
       +(\pa_3 u_3, \pa_3\eta v_i)_{G_1\cup G_{N,N+1}}-(\pa_3\eta u, \pa_3 {\bf a})_{G_1},
  \end{array}
\end{equation}
where we used that
$$(\na\eta u,\na{\bf a})_{G_1\cup G_{N,N+1}}=(\pa_3\eta u, \pa_3 {\bf a})=0$$
and
$$(\na' u_3,\na' v_i)_{G_{1,N}}=-(u_3,\Da' v_i)_{G_{1,N}}=-k_i\int_1^N\int_{\Si^i}u_3\,dx'dx_3=0.$$
Furthermore, we have
$$\begin{array}{rcl}
((u\cdot\na){\bf a},\eta {\bf a})
    &=&(\div( u{\bf a}),\eta{\bf a})\ek
             &=&-(u{\bf a},\na\eta{\bf a}+\eta\na{\bf a})=
                 -(u_3{\bf a},\pa_3\eta{\bf a})-(u'{\bf a},\eta\na'{\bf a})\ek
             &=&-(u_3v_i,\pa_3\eta v_i)_{G_1\cup G_{N,N+1}}-(\eta{\bf a},(u\cdot\na){\bf a})
\end{array}$$
 yielding
\begin{equation}
\label{E4.4} ((u\cdot\na){\bf a},\eta {\bf a})=-\ha(u_3v_i,\pa_3\eta
v_i)_{G_1\cup G_{N,N+1}}.
\end{equation}
Concerning the third term of the left-hand side of \eq{4.1}
we have
\begin{equation}
\label{E4.5} (({\bf a}\cdot\na)u,\eta {\bf a})=(v_i\pa_3u_3, \eta
v_i) =-(v_iu_3, \pa_3(\eta v_i))=-(u_3, \pa_3\eta v_i^2)_{G_1\cup
G_{N,N+1}}.
\end{equation}
The fourth term of \eq{4.1} is expanded as
\begin{equation}
\begin{array}{rcl}
\label{E4.6}
    ((u\cdot\na)u,\eta {\bf a})&=&((u^N\cdot\na)u,{\bf a})\ek
    &=&((u^N\cdot\na)u^N,{\bf a})-((u^N\cdot\na)[(\eta-1)u],{\bf a})_{G_1\cup
G_{N,N+1}},
\end{array}
\end{equation}
and finally for the pressure term we have
$$(p,\div(\eta{\bf a}))=(p,\pa_3\eta v_i).$$

Thus we get that
\begin{equation}
\label{E4.7}
    ((u^N\cdot\na)u^N,{\bf a}) = R+(p,\pa_3\eta v_i),
\end{equation}
where
\begin{equation}
\label{E4.8}
\begin{array}{l}
 R\equiv (f,\eta {\bf a})_{G_1}-(\na' u_3,\eta\na' v_i)_{G_1\cup G_{N,N+1}}
       -(\pa_3 u_3, \pa_3\eta v_i)_{G_1\cup G_{N,N+1}}\ek
       \hspace{1.5cm} + (\pa_3\eta u, \pa_3 {\bf a})_{G_1}
       +\ha(u_3v_i,\pa_3\eta v_i)_{G_1\cup G_{N,N+1}}\ek
       \hspace{1.5cm}+(u_3, \pa_3\eta v_i^2)_{G_1\cup G_{N,N+1}}
       +((u^N\cdot\na)[(\eta-1)u],{\bf a})_{G_1\cup G_{N,N+1}}.
\end{array}
\end{equation}
 In view of the construction of $\eta$ and \eq{3.30},
 we get from Corollary \ref{C3.2} that
 \begin{equation}
\label{E4.9}
\begin{array}{l}
 |R|\leq c(\Om,\Phi)
\end{array}
\end{equation}
with $c(\Om,\Phi)>0$ independent of $N$.

On the other hand, $u^N$ solves the system
\begin{equation}
\label{E4.10}
\begin{array}{rl}
  -\Da u^N  + (u^N\cdot\na){\bf a}+ ({\bf a}\cdot\na)u^N+(u^N\cdot\na)u^N+\na p =f_1 \,\,
  &\text{in }\Om,\ek
  \div u^N =g \,\, &\text{in }\Om,\ek
  u^N =0 \,\,&\text{on }\pa\Om,
  \end{array}
\end{equation}
in a weak sense, where
$$\begin{array}{l}
    f_1:= -\Da \eta u-2\na\eta\cdot\na u+({\bf a}\cdot\na)\eta u
        +(\eta^2-\eta)(u\cdot\na)u +\eta u_3\pa_3\eta u,\ek
    g:=\na\eta\cdot u=\pa_3\eta u_3.
\end{array}$$
From the construction of the cut-off function $\eta$, we get
    $$\supp f_1\subset G_1\cup G_{N,N+1},\quad \supp g\subset G_1\cup G_{N,N+1}$$
and, moreover,
$$\|f_1\|_{H^{-1}(\Om)}\leq c(\Om,\Phi)$$
in view of \eq{1.6} and Corollary \ref{C3.2}.
\par
Testing the first equation of \eq{4.10} with $u^N$, we get
\begin{equation}
\label{E4.11}
    \|\na u^N\|^2_{L^2(G_{N+1})}+((u^N\cdot\na){\bf a},u^N)
    +(({\bf a}\cdot\na)u^N,u^N)+((u^N\cdot\na)u^N, u^N)
    -(p,g)=(f_1,u^N).
\end{equation}
The second term in the left-hand side of \eq{4.11} is expanded as
\begin{equation}
\label{E4.12}
    ((u^N\cdot\na){\bf a},u^N)
    =-((u^N\cdot\na)u^N,{\bf a})-(g{\bf a},u^N),
\end{equation}
and the third term as
\begin{equation}
\label{E4.13}
    (({\bf a}\cdot\na)u^N,u^N)=(v_i\pa_3u^N,u^N)=
    -(v_i\pa_3u^N,u^N)=0
\end{equation}
since $v_i$ depends only on $x'$.
Concerning the fourth term, we get that
\begin{equation}
\label{E4.14}
    ((u^N\cdot\na)u^N, u^N)= -\ha(gu^N,u^N)
\end{equation}
 since $((u^N\cdot\na)u^N, u^N)= -((u^N\cdot\na)u^N,
u^N)-(gu^N,u^N)$.
Summarizing, we get from \eq{4.11}-\eq{4.14}
that
\begin{equation}
\label{E4.15}
\begin{array}{l}
    \|\na u^N\|^2_{L^2(G_{N+1})}-((u^N\cdot\na)u^N,{\bf a})\ek
    = (f_1,u^N)_{G_1\cup G_{N,N+1}}+(g{\bf
    a},u^N)_{G_1\cup G_{N,N+1}}+\ha(gu^N,u^N)_{G_1\cup G_{N,N+1}} + (p,g).
\end{array}
\end{equation}
Adding \eq{4.7} and \eq{4.15} yields
\begin{equation}
\label{E4.16}
\begin{array}{l}
    \|\na u^N\|^2_{L^2(G_{N+1})}= (f_1,u^N) + R\ek
  \hspace{2.5cm}   +(g{\bf a},u^N)_{G_1\cup G_{N,N+1}}+\ha(gu^N,u^N)_{G_1\cup G_{N,N+1}}+(p,\pa_3\eta v_i+g).
\end{array}
\end{equation}
Therefore, in view of the fact that
the constant in Poincar\'e's inequality for $G_{N+1}$
is independent of $N$ and depends only on the diameter of
$\Si^i$, we get that
$$\begin{array}{l}
    \|\na u^N\|^2_{L^2(G_{N+1})} \leq c(\Om)\|f_1\|_{H^{-1}(\Om)}\|\na u^N\|_{L^2(G_{N+1})} + R\ek
  \hspace{2.5cm}   +(g{\bf a},u^N)_{G_1\cup G_{N,N+1}}
  +\ha(gu^N,u^N)_{G_1\cup G_{N,N+1}}+(p,\pa_3\eta v_i+g).
\end{array}$$
Thus,
\begin{equation}
\label{E4.17}
  \|\na u^N\|^2_{L^2(G_{N+1})} \leq  |R| + |R_1| + (p,\pa_3\eta v_i+g)_{G_{N,N+1}},
\end{equation}
where
\begin{equation}
\label{E4.18}
    R_1 \equiv c(\Om)\|f_1\|^2_{H^{-1}(\Om)} +(g{\bf a},u^N)_{G_1\cup
    G_{N,N+1}}+\ha(gu^N,u^N)_{G_1\cup G_{N,N+1}}+(p,\pa_3\eta
    v_i+g)_{G_1}.
\end{equation}
Note that $|R_1| \leq c(\Om,\Phi)$ with $c(\Om,\Phi)$ independent of $N$,
which together with
\eq{4.9} yields
\begin{equation}
\label{E4.19}
  \|\na u^N\|^2_{L^2(G_{N+1})} \leq  C(\Om,\Phi) + (p,\pa_3\eta v_i+g)_{G_{N,N+1}}.
\end{equation}

Finally let us estimate $(p,\pa_3\eta v_i+g)_{G_{N,N+1}}$.
Let
$$\bar{p}(x):=p(x)-\frac{1}{|\Si^i|}\int_{\Si^i}p(x',x_3)\,dx'.$$
Then, using the Poincar\'e's inequality, Proposition \ref{P3.1}
and Corollary \ref{C3.2}, we have
$$\begin{array}{l}
    |(p,\pa_3\eta v_i+g)_{G_{N,N+1}}|\ek
     = \di\Big|\frac{1}{|\Si^i|}\int_N^{N+1}\big(\int_{\Si^i}p\,dx'\cdot\int_{\Si^i}\pa_3\eta
    (v_i+u_3)\,dx'\big)dx_3
    +\int_N^{N+1}(\bar{p},\pa_3\eta (v_i+u_3))_{\Si^i}\,dx_3\Big|\ek
    \leq\di\frac{|\Phi_i|}{|\Si^i|}\Big|\int_N^{N+1}\pa_3\eta\big(\int_{\Si^i}p\,dx'\big)dx_3\Big|
    +\int_N^{N+1}|\pa_3\eta||(\bar{p}, (v_i+u_3))_{\Si^i}|\,dx_3\ek
    =\di\frac{|\Phi_i|}{|\Si^i|}\Big|\int_{\Si^i}p(x',N)\,dx'+\int_{G_{N,N+1}}\eta\pa_3
    p\,dx\Big|
    +\int_N^{N+1}|\pa_3\eta||(\bar{p}, (v_i+u_3))_{\Si^i}|\,dx_3\ek
    \leq\di
    c(\Si^i,\Phi_i)\big|\int_{\Si^i}p(x',N)\,dx'\big|+c(\Om,\Phi)\ek
\hspace{5cm}    +\di c(\Si^i)\int_N^{N+1}|\pa_3\eta|\|\na'
p\|_{L^2(\Si^i)} \|v_i+u_3\|_{L^2(\Si^i)}\,dx_3\ek
    \leq \di c(\Om,\Phi)\Big(\big|\int_{\Si^i}p(x',N)\,dx'\big|+1\Big).
\end{array}$$
By the assumption \eq{4.0}, there is a subsequence $\{N_k\}\subset\N$
such that
    $$N_k\ra\infty (k\ra\infty),\quad
    \sup_{k\in\N}\big|\int_{\Si^i}p(x',N_k)\,dx'\big|<\infty.$$
Hence, we have
\begin{equation}
\label{E4.20}
    |(p,\pa_3\eta v_i+g)_{G_{N_k,N_k+1}}|<C(\Om,\Phi).
\end{equation}
Thus, we get finally from \eq{4.19} and \eq{4.20} that
$$    \|\na u^{N_k}\|^2_{L^2(\Om)} \leq C(\Om,\Phi),\quad\forall
    k\in\N,
$$
yielding \eq{4.2} in view of \eq{4.22}.
In particular, it follows from \eq{4.2} that
$\{u^{N}\}$ has a subsequence weakly convergent in $H^1_0(\Om)$
to some $\tilde{u}\in H^1_0(\Om)$ and hence $\|\tilde{u}\|_{H^1_0(\Om)}\leq C(\Om,\Phi)$.
By the way, due to \eq{4.22}, $\{u^{N}(x)\}$ converges to $u(x)$ for a.a. $x\in \Om^i\setminus G_1$,
which implies $u=\tilde{u}$ for a.a. $x\in\Om^i\setminus G_1$.
Thus we get
 $$u\in H^1_0(\Om), \quad\|u\|_{H^1_0(\Om)}\leq C(\Om,\Phi).$$
The proof of the theorem is complete.\hfill\qed

\begin{rem}
Let us obtain an equivalence condition for \eq{4.0}.
Consider the third equation w.r.t. $u_3$ in \eq{3.3}, that is,
$$-\Da u_3+(u\cdot \na){\bf a}_3+({\bf a}\cdot \na)u_3+(u\cdot \na)u_3+\pa_3p=f_3
$$
in $\Om^i$ for some $i\in\{1,\ldots,m\}$,
which may be simplified as (${\bf a}_3\equiv v$)
\begin{equation}
\label{E4.23}
-\Da u_3+\div(uv)+\pa_3(vu_3)+\div'(u'u_3)+\frac{1}{2}\pa_3(u_3^2)+\pa_3p=f_3.
\end{equation}
Now, integrating \eq{4.23} over $\Om^i_{1,N}$
 we have
\begin{equation}
\label{E4.24}
\int_{\Si^i} p(x',N)\,dx'-\int_{\Si^i} p(x',1)\,dx'
 =\int_1^N \int_{\pa\Si^i} \frac{\pa u_3}{\pa n'}\,dx'+R_2,
\end{equation}
where $n'$ denotes the unit outward normal vector at $\Si^i$ and $|R_2|\leq C(\Om,\Phi)$
in view of Proposition \ref{P3.1} and Corollary \ref{C3.2}.

 Therefore,
\begin{equation}
\label{E4.26}
\Big|\int_{\Si^i} p(x',N)\,dx'\Big|<C(\Om,\Phi) \quad\text{iff}\quad
\Big|\int_1^N \int_{\pa\Si^i} \frac{\pa u_3}{\pa n'}\,dx'\Big|<C(\Om,\Phi).
\end{equation}

\end{rem}


\begin{thebibliography}{99}


\bibitem{Am77} C.J. Amick, Steady solutions of the
Navier-Stokes equations in unbounded channels and pipes, Ann. Scuola
Norm. Sup. Pisa 4 (1977), 473-513

\bibitem{AmFr80} C.J. Amick and L.E. Fraenkel, Steady solutions of the Navier-Stokes
equations representing plane flows in channels of various types, Acta Math. 144  (1980), 81-152

\bibitem{Be04} H. Beir\~{a}o da Veiga, Time-periodic solutions of the Navier-Stokes
                equations in unbounded cylindrical domains:
                Leray's problem for periodic flows,
                Arch. Rational Mech. Anal. 178 (2005), 301-325

\bibitem{Bog79} M.E. Bogovskii, Solution of the first boundary
value problem for the equation of continuity of an incompressible
medium, Soviet Math. Dokl. 20 (1979), 1094-1098





\bibitem{Fa03} R. Farwig, Weighted $L_p$ Helmholtz decompositions
in infinite cylinders and in infinite layers. Adv. Diff. Equations 8 (2003), 357-384

\bibitem{FoFr00} M.A. Fontelos and A. Friedman, Stationary non-Newtonian fluid
flows in channel-like
and pipe-like domains. Arch. Ration. Mech. Anal. 151 (2000), 1-43

\bibitem{Ga94-1} G.P. Galdi, An introduction to the mathematical theory of the
Navier-Stokes equations, Vol. 1: Linearized steady problems,
Springer Tracts in Natural Philosophy, 38, Springer, 1994

\bibitem{Ga94-2} G.P. Galdi, An introduction to the mathematical theory of the
                  Navier-Stokes equations, Vol. II: Nonlinear steady problems,
                  Springer Tracts in Natural Philosophy, 39,  Springer, 1994


\bibitem{Ka85} L.P. Kapitanskii,Stationary solutions of the Navier-Stokes equations in periodic tubes.
(English Translation) J. Soviet Math. 28  (1985), 689-695

\bibitem{KaPi84}  L.P. Kapitanskii and K. Pileckas, On spaces of solenoidal vector fields and boundary
value problems for the Navier-Stokes equations in domains with noncompact boundaries.
(English Translation) Proc. Steklov Math. Inst. 159 (1984), 3-34

\bibitem{La59} O.A. Ladyzhenskaya, Stationary motion of a viscous incompressible fluid in a pipe.
Dokl. Akad. Nauk. SSSR 124, 551-553 (1959)

\bibitem{LS83} O.A. Ladyzhenskaya and V. A. Solonnikov,
Determination of solutions of boundary value problems for stationary
Stokes and Navier-Stokes equations having an unbounded Dirichlet
integral, J. Sov. Math. 21 (1983), 728-761

\bibitem{NaPi84} S.A. Nazarov and K. Pileckas, On the behavior of solutions of the Stokes
 and Navier-Stokes systems in domains with a periodically varying section. (English Translation)
Proc. Steklov Math. Inst. 159 (1984), 97-104

\bibitem{NaPi90} S.A. Nazarov and K. Pileckas, The Reynolds flow of a fluid in a three-dimensional
channel. (In Russian) Liet. Mat. Rink. 30  (1990), 772-783


\bibitem{Pi97} K. Pileckas,
          Strong solutions of the steady nonlinear
          Navier-Stokes system in domains with
          exits to infinity, Rend. Sem. Mat. Univ. Padova 97 (1997), 235-267

\bibitem{Pi00 } K. Pileckas, A. Sequeira and J.H. Videman, Steady flows of viscoelastic fluids in
domains with outlets to infinity. J. Math. Fluid Mech. 2 (2000), 185-218

\bibitem{RF07} M.-H. Ri and R. Farwig, Existence and
exponential stability in $L^r$-spaces of stationary Navier-Stokes
flows with prescribed flux in infinite cylindrical domains, Math.
Methods Appl. Sci. 30 (2007), 171-199


\end{thebibliography}
\end{document}